\documentclass{agtart_a}
\pdfoutput=1

\usepackage{xy,enumerate}
\xyoption{all}

\title{Rigidification of algebras over multi-sorted theories}

\author{Julia E Bergner}
\givenname{Julia E}
\surname{Bergner}
\address{Kansas State University\\138 Cardwell Hall\\\newline
Manhattan, KS 66506\\USA}
\email{bergnerj@member.ams.org}
\urladdr{}

\volumenumber{6}
\issuenumber{}
\publicationyear{2006}
\papernumber{68}
\startpage{1925}
\endpage{1955}

\doi{}
\MR{}
\Zbl{}


\received{9 August 2005}
\revised{8 September 2006}
\accepted{29 September 2006}
\published{14 November 2006}
\publishedonline{14 November 2006}
\proposed{}
\seconded{}
\corresponding{}
\editor{CPR}
\version{}

\arxivreference{math.AT/0508152}

\let\xysavmatrix\xymatrix
\def\xymatrix{\disablesubscriptcorrection\xysavmatrix}
\AtBeginDocument{\let\bar\wbar\let\tilde\wtilde\let\hat\what}

\makeatletter
\def\cnewtheorem#1[#2]#3{\newtheorem{#1}{#3}[section]
\expandafter\let\csname c@#1\endcsname\c@theorem}
\makeatother

\newenvironment{enumroman}{\begin{enumerate}[\upshape (i)]}
                                                {\end{enumerate}}
\theoremstyle{plain}
\newtheorem{theorem}{Theorem}[section]
\cnewtheorem{cor}[theorem]{Corollary}
\cnewtheorem{prop}[theorem]{Proposition}
\cnewtheorem{lemma}[theorem]{Lemma}

\theoremstyle{definition}
\cnewtheorem{definition}[theorem]{Definition}
\cnewtheorem{example}[theorem]{Example}
\newtheorem*{thank}{Acknowledgments}
\newtheorem*{notation}{Notation}

\newcommand{\Deltaop}{{\bf \Delta}^{op}}

\newcommand{\colim}{\mathrm{colim}}
\newcommand{\hocolim}{\mathrm{hocolim}}

\newcommand{\Hom}{\mathrm{Hom}}
\newcommand{\Map}{\mathrm{Map}}
\newcommand{\holim}{\mathrm{holim}}
\newcommand{\operad}{\mathrm{operad}}
\newcommand{\Algt}{\mathcal Alg^\mathcal T}
\newcommand{\SSetst}{\mathcal{SS}ets^\mathcal T}
\newcommand{\LSSetst}{\mathcal{LSS}ets^\mathcal T}
\newcommand{\SSetsd}{\mathcal{SS}ets^\mathcal D}
\newcommand{\alphau}{{\underline \alpha}}
\newcommand{\betau}{{\underline \beta}}

\newcommand{\SSets}{\mathcal {SS}ets}
\newcommand{\Sets}{\mathcal{S}ets}
\newcommand{\Tocat}{\mathcal T_{\mathcal {OC}at}}

\begin{document}

\begin{abstract}
We define the notion of a multi-sorted algebraic theory, which is
a generalization of an algebraic theory in which the objects are
of different ``sorts.'' We prove a rigidification result for
simplicial algebras over these theories, showing that there is a
Quillen equivalence between a model category structure on the
category of strict algebras over a multi-sorted theory and an
appropriate model category structure on the category of functors
from a multi-sorted theory to the category of simplicial sets. In
the latter model structure, the fibrant objects are homotopy
algebras over that theory.  Our two main examples of strict
algebras are operads in the category of simplicial sets and
simplicial categories with a given set of objects.
\end{abstract}

\maketitle

\section{Introduction}
Algebraic theories are useful in studying many standard algebraic
objects, such as monoids, abelian groups, and commutative rings.  An
algebraic theory provides a functorial means of describing particular
algebraic objects without specifying generating sets for the
operations to which the objects are subject, or for the relations
between these operations (Lawvere \cite{law}). Given a category
$\mathcal C$ of algebraic objects, the associated algebraic theory
$\mathcal T_{\mathcal C}$ (if it exists) is a small category with
products satisfying the property that specifying an object of
$\mathcal C$ is equivalent to giving a product-preserving functor
$\mathcal T_{\mathcal C} \rightarrow \Sets$.

Consider a category $\mathcal C$ with an associated algebraic
theory $\mathcal T$.  If a functor from $\mathcal T$ to the
category of simplicial sets preserves products, then it is
essentially a simplicial object in $\mathcal C$ and is thus a
combinatorial model for a topological object in $\mathcal C$, such
as a topological group when $\mathcal C$ is the category of
groups. We call such a functor a \emph{strict} $\mathcal
T$--\emph{algebra} (\fullref{strict}). If the functor
preserves products up to homotopy, we call it a \emph{homotopy}
$\mathcal T$--\emph{algebra} (\fullref{homotopy}). A
homotopy $\mathcal T$--algebra can be viewed as a simplicial set
with the appropriate algebraic structure ``up to homotopy," in a
higher-order sense. Using an appropriate notion of weak
equivalence on homotopy $\mathcal T$--algebras (Badzioch \cite[5.6]{bad}), the
following result relates strict and homotopy
$\mathcal T$--algebras:

\begin{theorem}[Badzioch {\cite[1.4]{bad}}] \label{rigid}
Let $\mathcal T$ be an algebraic theory.  Any homotopy $\mathcal
T$--algebra is weakly equivalent as a homotopy $\mathcal T$--algebra
to a strict $\mathcal T$--algebra.
\end{theorem}

As a motivation for the work in this paper, consider the category
of monoids.  There is an associated algebraic theory $\mathcal
T_M$, and thus a simplicial monoid can be specified by a $\mathcal
T_M$--algebra.  However, the notion of simplicial monoid can be
generalized to that of a simplicial category, by which we mean a
category enriched over simplicial sets, since a simplicial monoid
is a simplicial category with one object.  We would like to have a
generalization of Badzioch's theorem which applies to simplicial
categories. From the point of view of algebraic structure, the
main difference between a simplicial monoid and a simplicial
category with more than one object is that in the latter case the
description of the algebraic structure is more complicated, in
that two morphisms can be combined by the composition operation
only if they satisfy certain compatibility conditions on the
domain and range. Therefore, we would like to describe a more
general notion of theory which is capable of describing algebraic
structures in which the elements have various sorts or types, and
in which the operations which can be used to combine a collection
of elements depend on these sorts.

There is in fact such a ``multi-sorted" theory, $\Tocat$, such
that a product-preserving functor $\Tocat \rightarrow \Sets$ is
essentially a category with object set $\mathcal O$ (\fullref{ocat}). A simplicial category, analogously, can be viewed as
a product-preserving functor $\Tocat \rightarrow \SSets$.

A simpler example of an algebraic structure which requires the use
of a multi-sorted theory, which we will describe in more detail in
\fullref{group}, is that of a group acting on a set. There are
two sorts of elements, namely, the elements of the group and the
elements of the set.  Two elements of the group can be combined
via multiplication, or an element can be inverted. An element of
the group and an element of the set can be combined via the group
action.  However, the elements of the set cannot be combined with
one another in any nontrivial way, so the operations which we
allow depend on the sort of element involved. The example of a
module over a ring is constructed similarly in \fullref{ring}.

Another application of the notion of a multi-sorted theory gives a
convenient description of an operad.  In \fullref{operad}, we
characterize the theory $\mathcal T_{\operad}$ of operads.  An
operad in the category of sets is then a product-preserving
functor from $\mathcal T_{\operad}$ to the category of sets.  Thus,
we can describe an operad of sets as a diagram of sets given by
this multi-sorted theory.  We can similarly describe operads of
spaces.

A multi-sorted theory $\mathcal T$ is a category with products, so
we can define strict and homotopy $\mathcal T$--algebras as before
(see Definitions \ref{talg} and \ref{htalg}).  Using a definition
of weak equivalence for homotopy $\mathcal T$--algebras
(\fullref{wkequiv}), the main result which we prove for
multi-sorted theories is the following generalization of \fullref{rigid}:

\begin{theorem} \label{first}
Let $\mathcal T$ be a multi-sorted algebraic theory.  Any homotopy
$\mathcal T$--algebra is weakly equivalent as a homotopy $\mathcal
T$--algebra to a strict $\mathcal T$--algebra.
\end{theorem}

As Badzioch does, we actually prove a stronger statement in terms
of a Quillen equivalence of model category structures (\fullref{second}).

Using our example of the theory $\mathcal T_{\operad}$ of operads,
an operad in the category of simplicial sets is a strict $\mathcal
T_{\operad}$--algebra.  A homotopy operad, or sequence of simplicial
sets with the structure of an operad only up to homotopy, is then
a homotopy $\mathcal T_{\operad}$--algebra and can be rigidified to
a strict operad using this theorem.

Returning to the example of simplicial categories, let $\mathcal
O$ be a set and $\mathcal {SC}_{\mathcal O}$ the category of
simplicial categories with object set $\mathcal O$ in which the
morphisms are the identity on the objects.  In \cite{simpmon}, we
use \fullref{first} to prove a relationship between $\mathcal
{SC}_{\mathcal O}$ and the category of Segal categories with the
same set $\mathcal O$ in dimension zero.  In \cite{thesis}, we use
the ideas of this proof to prove an analogous relationship between
the category of all small simplicial categories and the category
of all Segal categories.

Throughout this paper, we frequently work in the category of
simplicial sets, $\SSets$.  Recall that a simplicial set is a
functor $\Deltaop \rightarrow \Sets$, where $\Delta$ denotes the
\emph{cosimplicial category} whose objects are the finite ordered
sets $[n]=(0, \ldots, n)$ and whose morphisms are the
order-preserving maps.  The \emph{simplicial category} $\Deltaop$
is then the opposite of this category.  Some examples of
simplicial sets are, for each $n \geq 0$, the $n$--simplex $\Delta
[n]$, its boundary $\dot \Delta [n]$, and, for any $0 \leq k \leq
n$, the simplicial set $V[n,k]$, which is $\dot \Delta [n]$ with
the $k$th face removed. More information about simplicial sets can
be found in Goerss and Jardine \cite[I.1]{gj}.

In this paper, we begin by recalling the definition of an
algebraic theory and stating some of its basic properties. Using
this definition as a model, we then define a multi-sorted theory.
We should note here that this notion is not a new one; similar
definitions are given by Ad\'{a}mek and Rosick\'{y}
\cite[3.14]{ar} and by Boardman and Vogt \cite[2.3]{bv}. (The
still more general definition of a finite limit theory is used by
Johnson and Walters \cite{jw}, and Rosick\'{y} proves a similar
result to \fullref{first} for limit theories \cite{ros}.)
Because our perspective is slightly different, however, we will
give a precise definition followed by some examples. Given a
multi-sorted theory $\mathcal T$, we define strict and homotopy
$\mathcal T$--algebras over a multi-sorted theory $\mathcal T$ and
show that the existence of a model category structure on the
category of all $\mathcal T$--algebras. We also show the existence
of a model category structure on the category of all functors
$\mathcal T \rightarrow \SSets$ in which the fibrant objects are
the homotopy $\mathcal T$--algebras. We then show that there is a
Quillen equivalence between these two model categories.

We note that the key here is the fact that we are considering
functors which preserve \emph{categorical} products.  An
interesting question which we hope to address further in future
work is whether such a rigidification holds in a category such as
chain complexes where the ``product" we are interested in, the
tensor product, is not a categorical product.

\begin{thank}
I am grateful to Bill Dwyer for suggesting this approach to
studying simplicial categories and operads.  I would also like to
thank Bernard Badzioch and Michael Johnson for helpful
conversations about this work, Ji\v{r}\'{i} Rosick\'{y} for
pointing out his related results, and the referee for numerous
suggestions for the improvement of this paper.  Partial support
from a Clare Boothe Luce Foundation Graduate Fellowship is also
gratefully acknowledged.
\end{thank}

\section{A summary of algebraic theories}

We first recall the definition of an ordinary algebraic theory.  More
details about algebraic theories can be found in Borceux \cite[Chapter
3]{bor}.

\begin{definition}
An \emph{algebraic theory} $\mathcal T$ is a small category with
finite products and which has as objects $T_n$ for $n \geq 0$
together with, for each $n$, an isomorphism $T_n \cong (T_1)^n$.
In particular, $T_0$ is the terminal object in $\mathcal T$.
\end{definition}

We can use theories to describe certain algebraic categories,
namely those which are determined by sets with $n$--ary operations
for each $n \geq 2$.  To do so, we need to use the notion of
adjoint pairs of functors.  Recall that a pair of functors
\[ \xymatrix@1{F\co  \mathcal C \ar@<.5ex>[r] & \mathcal D\,\,\colon\! R \ar@<.5ex>[l]}
\] is \emph{adjoint} (where $F$ is the left adjoint and $R$ is the right adjoint) if there is
a natural isomorphism
\[ \varphi_{X,Y}\co  \Hom_{\mathcal D}(FX,Y) \rightarrow
\Hom_{\mathcal C}(X,RY) \] for all objects $X$ in $\mathcal C$ and
$Y$ in $\mathcal D$.  The adjoint pair is sometimes written as the
triple $(F,R, \varphi)$ (Mac Lane \cite[IV.1]{macl}).

Now, consider a category $\mathcal C$ such that there exists a
forgetful functor
\[ \Phi \co \mathcal C \rightarrow Sets \]
taking an object of $\mathcal C$ to its underlying set, and its
left adjoint (a free functor)
\[ L\co Sets \rightarrow \mathcal C. \]  In other words, $\mathcal
C$ is required to have free objects.  If the category $\mathcal C$
and the adjoint pair $(\Phi, L)$ satisfy some additional technical
conditions (see \cite[3.9.1]{bor} for details), we will call
$\mathcal C$ an \emph{algebraic category}.

Given an object $X$ of an algebraic category $\mathcal C$, we have
a natural map
\[ \varepsilon_X\co  L\Phi(X) \rightarrow X \]
and given a set $A$, we have another map
\[ \eta_A\co A \rightarrow \Phi L(A). \]
In order to discuss a theory over the algebraic category $\mathcal
C$, consider a set $A$ together with a map $m_A\co  \Phi L(A)
\rightarrow A$ satisfying two conditions: the composite map
\[ \xymatrix@1{A \ar[r]^-{\eta_A} & \Phi L(A) \ar[r]^-{m_A} & A}
\]
is the identity map on $A$, and the diagram
\[ \xymatrix@1{(\Phi L)^2A \ar@<.5ex>[rr]^{\Phi L(m_A)} \ar@<-.5ex>[rr]_{\Phi \varepsilon_{LA}} && \Phi L(A)
\ar[r]^-{m_A} & A} \] is a coequalizer.  These maps define an
algebraic structure on the set $A$, specifically the structure
possessed by the objects of $\mathcal C$ \cite{law}.  (Note that
via this structure $\Phi L$ defines a monad on the category of
sets \cite[VI.1]{macl}.)

For example, if $\mathcal C = \mathcal G$, the category of groups,
$\Phi$ is the forgetful functor taking a group to its underlying
set, and $L$ is the free group functor taking a set to the free
group on that set, then these two conditions are precisely the
ones defining a group structure on the set $A$.

We would like to discuss the algebraic theory $\mathcal T$
corresponding to $\mathcal C$ to simplify this way of talking
about algebraic structure. Let $X$ be an object of $\mathcal C$.
We consider natural transformations of functors $\mathcal C
\rightarrow \Sets$
\[ \underbrace{\Phi (-) \times \cdots \times \Phi (-)}_n \rightarrow \Phi (-). \]
Using the adjointness of $\Phi$ and $L$, we have that
\[ \Phi (X) \cong \Hom_{\mathcal Sets}(\{1\}, \Phi (X)) \cong
\Hom_{\mathcal C}(L\{1\},X) \] where $\{1\}$ denotes the set with
one object, and we can think of $L\{1\}$ as the free object in
$\mathcal C$ on one generator, since $L$ is the free functor.
Hence, we have
\[ \begin{aligned}
\Phi (X)^n & = \Hom_{\mathcal Sets}(\{1\}, \Phi(X))^n \\
& = \Hom_{\mathcal Sets}(\coprod_n \{1\}, \Phi(X)) \\
& = \Hom_{\mathcal Sets}(\{1, \ldots ,n\}, \Phi (X)) \\
& = \Hom_{\mathcal C}(L\{1, \ldots ,n\},X).
\end{aligned} \]
Now, by Yoneda's Lemma we have a bijection between the set of
natural maps $\Phi (X)^n \rightarrow \Phi (X)$ and the set
$\Hom_{\mathcal C}(L\{1\},L\{1, \ldots ,n\})$. The objects
\[ L\{\phi\}=T_0, L\{1\}=T_1, \ldots ,L\{1, \ldots ,n\}=T_n, \ldots \]
are the objects of the algebraic theory $\mathcal T$ corresponding
to $\mathcal C$. The morphisms are the opposites of the ones in
$\mathcal C$ between these objects. More precisely stated,
$\mathcal T$ is the opposite of the full subcategory of
representatives of isomorphism classes of finitely generated free
objects of $\mathcal C$.
\vspace{3pt}

Given an object $X$ of $\mathcal C$, define a functor
$H_X\co \mathcal T \rightarrow \mathcal Sets$ such that
\[ H_X(L\{1, \ldots ,n\})= \Hom_{\mathcal C}(L\{1, \ldots ,n\},X) =
\Phi (X)^n. \] Now, the algebraic category $\mathcal C$ is
equivalent to the category of the functors $H_X$, namely, the full
subcategory of the category of functors $A\co \mathcal T \rightarrow
\Sets$ whose objects preserve products, or those for which the
canonical map $A(T_n) \rightarrow A(T_1)^n$ induced by the $n$
projection maps is an isomorphism of sets for all $n \geq 0$
\cite{law}.
\vspace{3pt}

\begin{example}
Let $\mathcal G$ denote the category of groups.  Consider the full
subcategory of $\mathcal G$ whose objects $T_n$ are the free
groups on $n$ generators for $n \geq 0$ (where $T_0$ is the
trivial group).  The opposite of this category is $\mathcal
T_{\mathcal G}$, the theory of groups.  It can be shown that the
category of product-preserving functors $\mathcal T_{\mathcal G}
\rightarrow \Sets$ is equivalent to the category $\mathcal G$.
\end{example}

\vspace{3pt}
Product-preserving functors from the theory $\mathcal T$ to
$\Sets$ are called \emph{algebras} over $\mathcal T$.  We would
also like to consider functors from an algebraic theory to the
category $\SSets$ of simplicial sets.  To do so, we must first
define a simplicial algebra over a theory $\mathcal T$.  For
simplicity, we will also use the term ``algebra" to refer to these
simplicial algebras.
\vspace{3pt}

\begin{definition}\cite[1.1]{bad}\label{strict}\qua
Given an algebraic theory $\mathcal T$, a \emph{(strict
simplicial)} $\mathcal T$--\emph{algebra} $A$ is a
product-preserving functor $A\co \mathcal T \rightarrow \mathcal
{SS}ets$.  Namely, the canonical map
\[ A(T_n) \rightarrow A(T_1)^n, \]
induced by the $n$ projection maps $T_n \rightarrow T_1$, is an
isomorphism of simplicial sets.  In particular, $A(T_0)$ is the
one-point space $\Delta [0]$.
\end{definition}

The category of all $\mathcal T$--algebras will be denoted $\Algt$.
Similarly, we have the notion of a homotopy algebra, for which we
only require products to be preserved up to homotopy:

\begin{definition}\cite[1.2]{bad}\label{homotopy}\qua
Given an algebraic theory $\mathcal T$, a \emph{homotopy}
$\mathcal T$--\emph{algebra} is a functor $X\co \mathcal T \rightarrow
\SSets$ which preserves products up to homotopy, ie, for each
$n$ the canonical map
\[ X(T_n) \rightarrow X(T_1)^n \]
is a weak equivalence of simplicial sets.  In particular, we
assume that $X(T_0)$ is weakly equivalent to $\Delta [0]$.
\end{definition}

There exists a forgetful functor, or evaluation map,
\[ U_{\mathcal T}\co  \Algt \rightarrow \mathcal {SS}ets \]
such that $U_{\mathcal T}(A)=A(T_1)$.  This functor has a left
adjoint, the free $\mathcal T$--algebra functor
\[ F_{\mathcal T}\co \mathcal{SS}ets \rightarrow \Algt \]
where, if $Y$ is any simplicial set,
\[ F_{\mathcal T}(Y)(T_1)=\coprod_{n \geq 0} \Hom_{\mathcal
T}(T_n,T_1) \times Y^n/\sim \] where the identifications come from
the structure of the algebraic theory \cite[2.1]{bad}.  More
specifically, if $\mathcal T_0$ denotes the initial theory (given
by representatives of isomorphism classes of finite sets), this
free functor is given by a coend
\[ F_{\mathcal T}(Y)(T_1) = \int^{\mathcal T_0} \Hom_{\mathcal T}(T_n,
T_1) \times Y^n \] as given by Schwede in \cite[2.3]{sch}.

\section{Multi-sorted algebraic theories}

We now generalize the definition of an algebraic theory to that of
a multi-sorted theory.

\begin{definition}
Given a set $S$, an $S$--\emph{sorted algebraic theory} (or
\emph{multi-sorted theory}) $\mathcal T$ is a small category with
objects $T_{\alphau^n}$ where $\alphau^n = <\alpha_1, \ldots
,\alpha_n>$ for $\alpha_i \in S$ and $n \geq 0$ varying, and such
that each $T_{\alphau^n}$ is equipped with an isomorphism
\[ T_{\alphau^n} \cong \prod_{i=1}^n T_{\alpha_i}. \]
For a particular $\alphau^n$, the entries $\alpha_i$ can repeat,
but they are not ordered.  In other words, $\alphau^n$ is a an
$n$--element subset with multiplicities.  There exists a terminal
object $T_0$ (corresponding to the empty subset of $S$).
\end{definition}
\vspace{2pt}

\begin{notation}
Lower-case Greek letters (with or without subscripts), say
$\alpha$ or $\alpha_i$, will be used to denote objects of $S$,
whereas underlined ones, say $\alphau^n$ or simply $\alphau$, will
denote an $n$--element subset of objects of $S$ (with
multiplicities) for $n \geq 1$.
\end{notation}
\vspace{2pt}

Notice that a theory with a single sort is a theory in the sense
of the previous section.
\vspace{2pt}

We would like to speak of multi-sorted theories corresponding to
categories which are analogous to the algebraic categories which
we had in the ordinary case. However, because we have several
objects (or ``sorts") $T_{\alpha}$ where we only had the object
$T_1$ in an ordinary theory, we have many pairs of adjoint
functors, one for each sort.
\vspace{2pt}

Let $\mathcal C$ a category with coproducts such that given any
element $\beta \in S$, we have a forgetful functor
\[ \Phi_\beta \co  \mathcal C \rightarrow \Sets \]
and its left adjoint, the free functor
\[ L_\beta \co  \Sets \rightarrow \mathcal C. \]

We would like the category $\mathcal C$ and these adjoint pairs to
satisfy the following analogous conditions to those of
\cite[3.9.1]{bor}:
\begin{enumerate}
\item The category $\mathcal C$ has coequalizers and kernel pairs
(ie, pullbacks of diagrams $X \rightarrow Y \leftarrow X$).
%

\item Each $\Phi_\beta$ reflects isomorphisms and preserves
regular epimorphisms (ie, those that are coequalizers).

\item For all $\beta \in S$, the composite functor $\Phi_\beta
L_\beta$ preserves filtered colimits.
\end{enumerate}
These conditions make $\mathcal C$ a kind of generalized algebraic
category.

Now, for each object $X$ in $\mathcal C$ and element $\beta \in
S$, we have a map
\[ \varepsilon_{X, \beta}\co L_{\beta}\Phi_{\beta}(X) \rightarrow X
\]
and, for each set $A$ a map
\[ \eta_{A, \beta}\co  A \rightarrow
\Phi_\beta L_\beta(A). \] As before, in order to make sense of the
notion of theory, we consider a set $A$ together with, for each
$\beta \in S$, a map
\[ m_{A, \beta}\co \Phi_\beta L_\beta(A) \rightarrow A \]
satisfying two conditions: the composite map
\[ \xymatrix@1{A \ar[r]^-{\eta_{A, \beta}} & \Phi_{\beta} L_{\beta}(A) \ar[r]^-{m_{A, \beta}} & A}
\]
is the identity map on $A$, and the diagram
\[ \xymatrix@1{(\Phi_\beta L_\beta)^2A \ar@<1ex>[rr]^{\Phi_\beta L_\beta(m_{A, \beta})}
\ar@<-1ex>[rr]_{\Phi_\beta \varepsilon_{L_\beta A, \beta}} &&
\Phi_\beta L_\beta (A) \ar[r]^-{m_{A, \beta}} & A} \] is a
coequalizer. These maps define a ``multi-sorted algebraic
structure" on $\mathcal C$.  In particular, we have a notion of
composition for certain elements of $\mathcal C$ depending on
their sorts.  Given this structure, we can now construct the
$S$--sorted theory corresponding to the category $\mathcal C$.

Given $\alpha_i, \beta \in S$, we consider natural transformations
of functors $\mathcal C \rightarrow \Sets$
\[ \Phi_{\alpha_1}(-) \times \cdots \times \Phi_{\alpha_n}(-)
\rightarrow \Phi_\beta (-). \] As before, we can apply these
functors to an object $X$ of $\mathcal C$ and rewrite to obtain a
map
\[ \Hom_{\mathcal Sets}(\{1\}, \Phi_{\alpha_1} (X)) \times \cdots
\times \Hom_{\mathcal Sets}(\{1\}, \Phi_{\alpha_n}(X)) \rightarrow
\Hom_{\mathcal Sets}(\{1\}, \Phi_\beta(X)) \] which, by
adjointness, is equivalent to
\[ \Hom_{\mathcal C}(L_{\alpha_1}\{1\},X) \times \cdots \times
\Hom_{\mathcal C}(L_{\alpha_n}\{1\}, X) \rightarrow \Hom_{\mathcal
C}(L_\beta\{1\},X). \] Since ${\mathcal C}$ has coproducts, we can
rewrite this map as
\[ \Hom_{\mathcal C}(L_{\alpha_1}\{1\} \amalg \cdots \amalg
L_{\alpha_n}\{1\},X) \rightarrow \Hom_{\mathcal C}(L_\beta\{1\},X).
\]  Then, by Yoneda's Lemma, there is a bijection between the set
of natural transformations
\[ \Phi_{\alpha_1}(-) \times \cdots \times \Phi_{\alpha_n}(-)
\rightarrow \Phi_\beta(-) \] and the set
\[ \Hom_{\mathcal C}(L_\beta\{1\}, \coprod_{k=1}^n L_{\alpha_k}\{1\}). \]  The objects of the
theory $\mathcal T$ corresponding to ${\mathcal C}$ are given by
finite coproducts of ``free" objects $L_{\alpha_k}\{1\}$ of
${\mathcal C}$ for all choices of $\alpha_k$, and the morphisms are
the opposites of those of ${\mathcal C}$. Let $X$ be an object of
${\mathcal C}$ and $(\alpha_1, \ldots ,\alpha_n) \in S^n$ an
$n$--tuple of elements in $S$. We define the functor $H_{X,
\alpha_1, \ldots , \alpha_n}\co  \mathcal T \rightarrow \mathcal
Sets$ such that
\[ H_{X, \alpha_1, \ldots, \alpha_n} (\coprod_{k=1}^n L_{\alpha_k}\{1\}) = \Hom_{\mathcal C}(\coprod_{k=1}^n L_{\alpha_k}\{1\}, X)
= \Phi_{\alpha_1}(X) \times \cdots \times \Phi_{\alpha_n}(X). \]
(Note that we still write the ``coproduct" to denote an object of
$\mathcal T$ to be consistent with previous notation, even though
in $\mathcal T$ it is actually a product.) The category $\mathcal
C$ is equivalent to the category of all such functors if it
satisfies the conditions given above.

We now consider some examples.

\begin{example} \label{group}
Consider pairs $(G,X)$, where $G$ is a group and $X$ is a set. We
can obtain two different 2--sorted theories from these pairs, one
corresponding to the category of unstructured pairs, and the other
corresponding to the category of pairs $(G,X)$ with a given action
of the group $G$ on the set $X$.

In each case, we have two forgetful functors and their respective
left adjoints. We begin with the category of unstructured pairs,
which we denote $\mathcal P$.  The objects are the pairs $(G,X)$
and the morphisms $(G,X) \rightarrow (H,Y)$ consist of pairs
$(\varphi, f)$ where $\varphi\co G \rightarrow H$ is a group
homomorphism and $f\co X \rightarrow Y$ is a map of sets. For each
sort $i=1,2$ we have a forgetful map
\[ \Phi_i\co  \mathcal P \rightarrow \Sets \]
and its left adjoint
\[ L_i\co  \Sets \rightarrow \mathcal P. \]
When $i=1$, we have, for any group $G$ and set $X$,
\[ \Phi_1(G,X) = G \]
(where on the right-hand side $G$ denotes the underlying set of
the group $G$) and for any set $S$
\[ L_1(S)= (F_S, \phi) \]
where $F_S$ denotes the free group on the set $S$.

Similarly, when $i=2$, we define
\[ \Phi_2(G,X) = X \]
and
\[ L_2(S) = (e, S) \]
where $e$ denotes the trivial group.

In order to determine the objects of our theory, consider functors
\[ F_{i,j}\co  \mathcal P \rightarrow \mathcal Sets \] such that
$F_{i,j}(G,X) = G^i \times X^j$.  In other words,
\[ F_{i,j}(G,X) = \Hom_{\mathcal P}(L_1\{1, \ldots ,i\} \amalg L_2\{1,
\ldots ,j\}, (G,X)) \] where $\{1, \cdots ,i\}$ denotes the set
with $i$ elements and similarly for $\{1, \ldots ,j\}$. The
objects of the theory will be representatives of the isomorphism
classes of the $L_1\{1, \ldots, i\}$ $ \amalg L_2\{1, \ldots, j\}$
for all choices of $i$ and $j$. This coproduct in $\mathcal P$ is
defined to be the coproduct of each element in the pairs. Thus we
have
\[ (G,X) \amalg (G',X') = (G \ast G', X \amalg X') \] where $G \ast
G'$ denotes the free product of groups.  So, our corresponding
theory is the opposite of the full subcategory of $\mathcal P$
whose objects are of the form $L_1\{1, \ldots, i\} \amalg L_2\{1,
\ldots, j\}$.

When we equip each pair $(G,X)$ with an action of $G$ on $X$ to
obtain another category which we denote $\mathcal{PA}$, the
process is identical until we have to specify the coproduct, since
in this case we need to take the group actions into account. We
then have the coproduct in $\mathcal{PA}$
\[ (G,X) \amalg (G',X') = (H, (H \times_G X) \amalg
(H \times_{G'} X')) \] where $H= G \ast G'$ and we have defined
\[ H \times_G X = \{(h,x)|h \in H, x \in X\}/\sim  \]
when $(hg,x) \sim (h,gx)$ for any $g \in G$.  We can now take the
opposite of a full subcategory of $\mathcal{PA}$ as above to
obtain the corresponding theory.  In particular, the objects of
the theory look like
\[ L_1\{1, \ldots ,i\} \amalg L_2\{1, \ldots, j\}
= (F_i, F_i \times \{1, \ldots ,j\}), \] where $F_i$ denotes the
free group on $i$ generators.
\end{example}

\begin{example} \label{ring}
A very similar example is the case of a commutative ring $R$ and
an $R$--module $A$.  Again, we have two different 2--sorted
theories: one where we have a ring $R$ and regard $A$ merely as an
abelian group, and the other where we consider the $R$--module
structure on $A$.

As before, we begin with $\mathcal{PR}$, the category of pairs
with no additional structure.  We have the forgetful map
\[ \Phi_1\co  \mathcal{PR} \rightarrow \mathcal Sets \]
where $\Phi_1(R,A) = R$ for any ring $R$ and abelian group $A$,
where on the right side $R$ is the underlying set of the ring $R$.
Its left adjoint is the functor
\[ L_1\co  \mathcal Sets \rightarrow \mathcal{PR} \]
where for any set $S$, $L_1(S) = (\mathbb Z[S], e)$, where
$\mathbb Z[S]$ is the free commutative ring on the set $S$ and $e$
denotes the trivial (abelian) group.  Then we have the map
\[ \Phi_2\co  \mathcal{PR} \rightarrow \mathcal Sets \]
such that $\Phi_2(R,A) = A$, where again on the right hand side
$A$ is the underlying set of the abelian group $A$.  Its left
adjoint is the map
\[ L_2\co  \mathcal Sets \rightarrow \mathcal{PR} \]
where $L_2(S) = (\mathbb Z, FA_S)$ where $FA_S$ denotes the free
abelian group on the set $S$.

To know what the objects of this 2--sorted theory are, we need to
know what the coproduct is.  We have that
\[ (R,A) \amalg (R',A') = (R \otimes_{\mathbb Z} R', A \oplus A'),
\] and from there we can obtain a theory as in the previous
example.

Now consider the category $\mathcal{PM}$ whose objects are pairs
$(R,A)$ where $R$ is a ring and $A$ is a module over $A$. If $A$
and $A'$ are modules over $R$ and $R'$, respectively, we have a
coproduct similar to that in the group action example. So, we say
that
\[ (R,A) \amalg (R',A') = (R \otimes_{\mathbb Z} R',(R' \otimes_{\mathbb Z} A) \oplus (R \otimes_{\mathbb Z} A')) \]
and construct the corresponding theory as before.
\end{example}

\begin{example} \label{operad}
Another example of a multi-sorted theory is the $\mathbb N$--sorted
theory of symmetric operads.  Recall that an operad in the
category of sets is a sequence of sets $\{P(k)\}_{k \geq 0}$, a
unit element $1 \in P(1)$, each with a right action of the
symmetric group $\Sigma_k$, and operations
\[ P(k) \times P(j_1) \times \cdots \times P(j_k) \rightarrow P(j_1+ \cdots +j_k) \] satisfying
associativity, unit, and equivariance conditions
\cite[II.1.4]{mss}.

There is a notion of a free operad on $n$ generators at levels $m_1,
\ldots ,m_n$ (Markl, Shnider and Stasheff \cite[Section II.1.9]{mss},
Rezk \cite[2.3.6]{rezkop}).  Specifically, such a free operad has, for each
$1 \leq i \leq n$, a generator in $P(m_i)$.  Note that the values of
$m_i$ can repeat. For example, one can think of the free operad on $n$
generators, each at level 1, as the free monoid on $n$ generators.

In the category of operads, consider the full subcategory of
isomorphism classes of free operads.  Each object in this
category, then, can be described as the free operad on $n$
generators at levels $m_1, \ldots ,m_n$ for some $n \geq 0$ and
$m_1, \ldots ,m_n$. The opposite of this category is the theory of
operads. Using the notation we have set up for multi-sorted
theories, we have that $T_\alpha$ for $\alpha \in \mathbb N$ is
just the free operad on one generator at level $\alpha$ and for
$\alphau^n = <\alpha_1, \ldots ,\alpha_n>$, we have that
$T_{\alphau^n}$ is the free operad on $n$ generators at levels
$\alpha_1, \ldots \alpha_n$.

There is also a notion of non-symmetric (or non--$\Sigma$) operads,
where we no longer have an action of the symmetric group or an
equivariance condition \cite[II.1.14]{mss}.  We can define the
theory of non--$\Sigma$ operads analogously, taking the opposite of
the full subcategory of isomorphism classes of free non--$\Sigma$
operads in the category of all non--$\Sigma$ operads.
\end{example}

\begin{example} \label{ocat}
Consider the category $\mathcal {OC}at$ whose objects are the
categories with a fixed object set $\mathcal O$ and whose
morphisms are the functors which are the identity map on the
objects. There is a theory $\mathcal T_{\mathcal {OC}at}$
associated to this category. The objects of the theory are
isomorphism classes of categories which are freely generated by
directed graphs with vertices corresponding to the elements of the
set $\mathcal O$. This theory will be sorted by pairs of elements
in $\mathcal O$, corresponding to the morphisms with source the
first element and target the second.  In other words, this theory
is $(\mathcal O \times \mathcal O)$--sorted.

In particular, consider $\alpha = (x,y)\in \mathcal O \times
\mathcal O$.  Then, if $x \neq y$, $T_\alpha$ is the category with
object set $\mathcal O$ and one nonidentity morphism with source
$x$ and target $y$.  If $x=y$, then $T_\alpha$ is the category
freely generated by one morphism from $x$ to itself and no other
nonidentity morphisms.

In general, if $\alphau =<\alpha_1, \ldots ,\alpha_n>$, then
$T_\alphau$ is the category with object set $\mathcal O$ and
morphisms freely generated by the morphisms given for each
$\alpha_i$ as in the previous case.

Consider the forgetful functor $\Phi_{\alpha}\co  \mathcal {OC}at
\rightarrow \mathcal Sets$ where, for any object $X$ in $\mathcal
{OC}at$ and $\alpha=(x,y)$,
\[ \Phi_{\alpha}(X)= \Hom_X(x,y). \]
Its left adjoint then is the free functor $L_\alpha$ defined by,
for a set $A$,
\[ \begin{aligned}
L_\alpha (A)=
\begin{cases}
C \text{ with } \Hom_C(x,y) = A & \text{if } x \neq y \\
C \text{ with } \Hom_C (x,y) = F_A & \text{if } x =y
\end{cases}
\end{aligned} \]
where $F_A$ is the free monoid generated by the set $A$ and where
in each case there are no other nonidentity morphisms in the
category $C$.
\end{example}
\vspace{-2pt}

As with ordinary algebraic theories, we can define strict and
homotopy $\mathcal T$--algebras for a multi-sorted theory $\mathcal
T$.
\vspace{-2pt}

\begin{definition} \label{talg}
Given an $S$--sorted theory $\mathcal T$, a \emph{(strict
simplicial)} $\mathcal T$--\emph{algebra} is a product-preserving
functor $A\co \mathcal T \rightarrow \mathcal{SS}ets$.  Here,
product-preserving means that the canonical map
\[ A(T_{\alphau^n}) \rightarrow \prod_{i=1}^n A(T_{\alpha_i}), \]
induced by the projections $T_{\alphau^n} \rightarrow
T_{\alpha_i}$ for all $1 \leq i \leq n$, is an isomorphism of
simplicial sets.
\end{definition}
\vspace{-2pt}

As before, we will denote the category of strict $\mathcal
T$--algebras by $\Algt$.
\vspace{-2pt}

\begin{definition} \label{htalg}
Given an $S$--sorted theory $\mathcal T$, a \emph{homotopy}
$\mathcal T$--\emph{algebra} is a functor $X\co \mathcal T \rightarrow
\mathcal{SS}ets$ which preserves products up to homotopy, ie,
for all $\alpha \in S^n$, the canonical map
\[ X(T_{\alphau^n}) \rightarrow \prod_{i=1}^n X(T_{\alpha_i}) \]
induced by the projection maps $T_{\alphau^n} \rightarrow
T_{\alpha_i}$ (for each $1 \leq i \leq n$) is a weak equivalence
of simplicial sets.
\end{definition}
\vspace{-2pt}

We would like to prove a rigidification result similar to \fullref{rigid} above.  We begin by finding model category structures
for $\mathcal T$--algebras and homotopy $\mathcal T$--algebras.  We
then find a Quillen equivalence between these model category
structures $\mathcal T$--algebras for any multi-sorted theory
$\mathcal T$.
\vspace{-2pt}

\section{Model category structures}

In this section, we define, given a multi-sorted theory $\mathcal
T$, model category structures on the category of diagrams
$\mathcal T \rightarrow \SSets$ and on the category of $\mathcal
T$--algebras. We begin with a review of model category structures.

Recall that a model category structure on a category ${\mathcal C}$
is a choice of three distinguished classes of morphisms:
fibrations, cofibrations, and weak equivalences.  A (co)fibration
which is also a weak equivalence will be called an \emph{acyclic
(co)fibration}.  With this choice of three classes of morphisms,
${\mathcal C}$ is required to satisfy the following five axioms
(Dwyer and Spali{\'n}ski \cite[3.3]{ds}).
\begin{itemize}\leftskip 10pt
\item[(MC1)] ${\mathcal C}$ has all small limits and colimits.

\item[(MC2)] If $f$ and $g$ are maps in ${\mathcal C}$ such that
their composite $gf$ exists, then if two of $f$, $g$, and $gf$ are
weak equivalences, then so is the third.

\item[(MC3)] If a map $f$ is a retract of $g$ and $g$ is a
fibration, cofibration, or weak equivalence, then so is $f$.

\item[(MC4)] If $i\co A \rightarrow B$ is a cofibration and $p\co X
\rightarrow Y$ is a fibration, then a dotted arrow lift exists in
any solid arrow diagram of the form
\[ \xymatrix{A \ar[r] \ar[d]^i & X \ar[d]^p \\
B \ar[r] \ar@{-->}[ur] & Y} \] if either
\begin{enumroman}
\item $p$ is a weak equivalence, or

\item $i$ is a weak equivalence.
\end{enumroman}
(In this case we say that $i$ has the \emph{left lifting property}
with respect to $p$ and that $p$ has the \emph{right lifting
property} with respect to $i$.)

\item[(MC5)] Any map $f$ can be factored two ways:
\begin{enumroman}
\item $f= pi$ where $i$ is a cofibration and $p$ is an acyclic
cofibration, and

\item $f=qj$ where $j$ is an acyclic cofibration and $p$ is a
fibration.
\end{enumroman}
\end{itemize}

An object $X$ in ${\mathcal C}$ is \emph{fibrant} if the unique map
$X \rightarrow \ast$ from $X$ to the terminal object is a
fibration.  Dually, $X$ is \emph{cofibrant} if the unique map
$\phi \rightarrow X$ from the initial object to $X$ is a
cofibration.  The factorization axiom MC5 guarantees that each
object $X$ has a weakly equivalent fibrant replacement $\widehat
X$ and a weakly equivalent cofibrant replacement $\widetilde X$.
These replacements are not necessarily unique, but they can be
chosen to be functorial in the cases we will use
Hovey \cite[1.1.3]{hovey}.

The model category structures which we will discuss are all
cofibrantly generated.  In a cofibrantly generated model category,
there are two sets of morphisms, one of generating cofibrations and
one of generating acyclic cofibrations, such that a map is a fibration
if and only if it has the right lifting property with respect to the
generating acyclic cofibrations, and a map is an acyclic fibration if
and only if it has the right lifting property with respect to the
generating cofibrations (Hirschhorn \cite[11.1.2]{hirsch}).  To
describe such model categories, we make the following definition.

We are now able to state the theorem, due to Kan, that we will use
to prove our model category structures in this paper.

\begin{theorem}{\rm\cite[11.3.2]{hirsch}}\label{CofGen}\qua
Let $\mathcal M$ be a cofibrantly generated model category with
generating cofibrations $I$ and generating acyclic cofibrations
$J$.  Let $\mathcal N$ be a category that satisfies axiom MC1 such
that there exists a pair of adjoint functors
\[ F\co  \mathcal M
\rightleftarrows \mathcal N\,\,\colon\! U. \]  If $FI=\{Fu \mid u \in I\}$ and
$FJ=\{Fv \mid v \in J\}$, and if
\begin{enumerate}[\upshape(1)]
\item each of the sets $FI$ and $FJ$ permits the small object
argument \cite[10.5.15]{hirsch}, and

\item $U$ takes (possibly transfinite) colimits of pushouts along
maps in $FJ$ to weak equivalences,
\end{enumerate}
then there is a cofibrantly generated model category structure on
$\mathcal N$ for which $FI$ is a set of generating cofibrations
and $FJ$ is a set of generating acyclic cofibrations, and the weak
equivalences are the maps that $U$ sends to weak equivalences in
$\mathcal M$.
\end{theorem}

We will refer to the standard model category structure on the
category $\SSets$ of simplicial sets.  In this case, a weak
equivalence is a map of simplicial sets $f\co X \rightarrow Y$ such
that the induced map $|f|\co|X| \rightarrow |Y|$ is a weak homotopy
equivalence of topological spaces.  The cofibrations are
monomorphisms, and the fibrations are the maps with the right
lifting property with respect to the acyclic cofibrations
\cite[I.11.3]{gj}. This model category structure is cofibrantly
generated; a set of generating cofibrations is $I=\{\dot \Delta
[n] \rightarrow \Delta [n] \mid n \geq 0\}$, and a set of
generating acyclic cofibrations is $J=\{V[n,k] \rightarrow \Delta
[n] \mid n \geq 1, 0 \leq k \leq n\}$.

We will also need the notion of a simplicial model category
$\mathcal M$. For any objects $X$ and $Y$ in a simplicial category
$\mathcal M$, the \emph{function complex} is the simplicial set
$\Map(X,Y)$.

\begin{definition}\cite[9.1.6]{hirsch}\qua
A \emph{simplicial model category} $\mathcal M$ is a model
category $\mathcal M$ that is also a simplicial category such that
the following two axioms hold:
\begin{itemize}\leftskip 10pt
\item[(SM6)] For every two objects $X$ and $Y$ of $\mathcal M$ and
every simplicial set $K$, there are objects $X \otimes K$ and
$Y^K$ in $\mathcal M$ such that there are isomorphisms of
simplicial sets
\[ \Map(X \otimes K, Y) \cong \Map(K, \Map(X,Y)) \cong \Map(X,Y^K)
\] that are natural in $X$, $Y$, and $K$.

\item[(SM7)] If $i\co A \rightarrow B$ is a cofibration in $\mathcal
M$ and $p\co X \rightarrow Y$ is a fibration in $\mathcal M$, then
the map of simplicial sets
\[ i^*\times p_*\co\Map(B,X) \rightarrow \Map(A,X)
\times_{\Map(A,Y)} \Map(B,Y) \] is a fibration which is an acyclic
fibration if either $i$ or $p$ is a weak equivalence.
\end{itemize}
\end{definition}
\vspace{-2pt}

It is important to note that a function complex in a simplicial
model category is only homotopy invariant in the case that $X$ is
cofibrant and $Y$ is fibrant. For the general case, we have the
following definition:
\vspace{-2pt}

\begin{definition}\cite[17.3.1]{hirsch}\qua
A \emph{homotopy function complex} $\Map^h(X,Y)$ in a simplicial
model category $\mathcal M$ is the simplicial set $\Map(\widetilde
X, \widehat Y)$ where $\widetilde X$ is a cofibrant replacement of
$X$ in $\mathcal M$ and $\widehat Y$ is a fibrant replacement for
$Y$.
\end{definition}
\vspace{-2pt}

Several of the model category structures that we use are obtained
by localizing a given model category structure with respect to a
map or a set of maps.  Suppose that $P = \{f\co A \rightarrow B\}$ is
a set of maps with respect to which we would like to localize a
model category $\mathcal M$.
\vspace{-2pt}

\begin{definition} \label{local}
A $P$--\emph{local} object $W$ is a fibrant object of $\mathcal M$
such that for any $f\co A \rightarrow B$ in $P$, the induced map on
homotopy function complexes
\[ f^*\co\Map^h(B,W) \rightarrow \Map^h(A,W) \]
is a weak equivalence of simplicial sets.  A map $g\co X \rightarrow
Y$ in $\mathcal M$ is then a $P$--\emph{local equivalence} if for
every local object $W$, the induced map on homotopy function
complexes
\[ g^*\co  \Map^h(Y,W) \rightarrow \Map^h(X,W) \]
is a weak equivalence of simplicial sets.
\end{definition}

Given a multi-sorted theory $\mathcal T$, let
$\mathcal{SS}ets^\mathcal T$ denote the category of functors
$\mathcal T \rightarrow \mathcal{SS}ets$. Note that the category
$\Algt$ of strict $\mathcal T$--algebras is a full subcategory of
$\mathcal{SS}ets^\mathcal T$.

The category $\SSetst$ is an example of a category of diagrams. In
general, given any small category $\mathcal D$, there is a
category $\SSetsd$ of $\mathcal D$--diagrams in $\SSets$, or
functors $\mathcal D \rightarrow \SSets$.  We can obtain two model
category structures on $\SSetsd$ by the following results.

\begin{theorem}{\rm \cite[IX 1.4]{gj}}\qua
Given the category $\mathcal{SS}ets^\mathcal D$ of $\mathcal
D$--diagrams of simplicial sets, there is a simplicial model
category structure $\mathcal{SS}ets^\mathcal D_f$ in which the
weak equivalences and fibrations are objectwise and in which the
cofibrations are the maps which have the left lifting property
with respect to the maps which are both fibrations and weak
equivalences.
\end{theorem}

\begin{theorem}{\rm  \cite[VIII 2.4]{gj}}\qua
There is a simplicial model category $\mathcal{SS}ets^\mathcal
D_c$ in which the weak equivalences and the cofibrations are
objectwise and in which the fibrations are the maps which have the
right lifting property with respect to the maps which are
cofibrations and weak equivalences.
\end{theorem}

We now return to the situation where our small category is a
multi-sorted theory $\mathcal T$.  We would like to have an
evaluation map and its left adjoint as in the ordinary case (see
the end of section 2 above), but here we will have one for each
$\alpha \in S$. These evaluation maps look like
\[ U_\alpha\co  \Algt \rightarrow \mathcal{SS}ets \]
such that
\[U_\alpha (A)=A(T_\alpha) \]
for any $\mathcal T$--algebra $A$.

Each functor $U_\alpha$ has a left adjoint, the free functor
\[ F_\alpha \co \mathcal {SS}ets \rightarrow \Algt \]
such that, given a simplicial set $Y$ and object $T_\betau$ in
$\mathcal T$,
\[ F_{\alpha}(Y)(T_\betau)= \coprod_{n \geq 0} (\Hom_{\mathcal
T}(T_{\alpha, \ldots ,\alpha}, T_\betau) \times Y^n)/\sim. \]  As
before, this free functor can be defined precisely as a coend over
the initial (single-sorted) theory, regarded as the subcategory of
the initial $S$--sorted theory whose objects are $(T_\alpha)^n$ for
$n \geq 0$,
\[ F_{\alpha}(Y)(T_\betau)= \int^{\mathcal T_0} \Hom_{\mathcal
T}((T_\alpha)^n, T_\betau) \times Y^n. \]

Given a theory $\mathcal T$ (regular or multi-sorted), define a
weak equivalence in the category $\Algt$ of $\mathcal T$--algebras
to be a map which induces a weak equivalence of simplicial sets
after applying the evaluation functor $U_{\alpha}$ for each sort
$\alpha$. Similarly, define a fibration of $\mathcal T$--algebras
to be a map $f$ such that $U_{\alpha}(f)$ is a fibration of
simplicial sets for all $\alpha$. Then define a cofibration to be
a map with the left lifting property with respect to the maps
which are fibrations and weak equivalences.

The following theorem is a generalization of a result by Quillen
\cite[II.4]{quillen}.

\begin{theorem}
Let $\mathcal T$ be an $S$--sorted theory.  There is a cofibrantly
generated model category structure on $\Algt$ with the weak
equivalences, fibrations, and cofibrations as defined above.
\end{theorem}

\begin{proof}
We use a slightly generalized version of \fullref{CofGen} with
the adjoint pairs $F_\alpha\co  \SSets \leftrightarrows \Algt\,\,\colon\!
U_\alpha$ for all $\alpha \in S$ and using the cofibrantly
generated model structure on $\SSets$ as given above. The
existence of limits and colimits follows just as they do in the
case where $\mathcal T$ is an ordinary theory
\cite[II.4]{quillen}.  Thus, verifying conditions (1) and (2) will
result in a model structure on $\Algt$ for which the set
$FI=\{F_\alpha \dot \Delta [n] \rightarrow F_\alpha \Delta [n]
\mid \alpha \in S, n \geq 0 \}$ is a set of generating
cofibrations and $FJ=\{F_\alpha V[n,k] \rightarrow F_\alpha \Delta
[n] \mid \alpha \in S, n \geq 1, 0 \leq k \leq n\}$ is a set of
generating acyclic cofibrations.

We first show that $FI$ and $FJ$ satisfy the small object
argument. Consider some $\mathcal T$--algebra $A$, which can be
written as a directed colimit $\colim_m (A_m)$ and can therefore
be computed objectwise. Thus, we can show that $F_\alpha \dot
\Delta [n]$ is small:
\[ \begin{aligned}
\Hom_{\Algt}(F_\alpha \dot \Delta [n], \colim_m (A_m)) & =
\Hom_{\SSets}(\dot \Delta [n], U_\alpha \colim_m (A_m)) \\
& = \Hom_{\SSets}(\dot \Delta [n], \colim_m (U_\alpha A_m)) \\
& = \colim_m \Hom_{\SSets}(\dot \Delta [n], U_\alpha A_m) \\
& = \colim_m \Hom_{\Algt}(F_\alpha \dot \Delta [n], A_m).
\end{aligned} \]
The object $V[n,k]$ can be shown to be small analogously, so we
have proved statement (1).

To prove statement (2), we need to show that taking a pushout
along a map in $FJ$ results in a map which is a weak equivalence
in $\Algt$.  Note that since weak equivalences are taken
levelwise, a (transfinite) directed colimit of weak equivalences
is still a weak equivalence, so checking a single pushout
suffices.

Consider a map $F_\alpha V[n,k] \rightarrow F_\alpha \Delta [n]$
in $FJ$ and a map $F_\alpha V[n,k] \rightarrow A$ for some object
$A$ of $\Algt$.  We then take the pushout $B$ in the following
diagram:
\[ \xymatrix{F_\alpha V[n,k] \ar[r] \ar[d] & A \ar[d] \\
F_\alpha \Delta [n] \ar[r] & B}. \] Suppose that $X \rightarrow Y$
is a map in $\Algt$ with the right lifting property with respect
to the maps in $FJ$.  Note by adjointness that it is just a
levelwise fibration of simplicial sets.  Then in the diagram
\[ \xymatrix{F_\alpha V[n,k] \ar[r] \ar[d] & A \ar[d] \ar[r] & X \ar[d] \\
F_\alpha \Delta [n] \ar[r] & B \ar[r] & Y} \] a lift $F_\alpha
\Delta [n] \rightarrow X$ exists, which implies by universality
that there is also a lift $B \rightarrow X$.

Now consider the diagram
\[ \xymatrix{V[n,k] \ar[r] \ar[d] & U_\alpha A \ar[d] \ar[r] & U_\alpha X \ar[d] \\
\Delta [n] \ar[r] & U_\alpha B \ar[r] & U_\alpha Y} \] where the
left-hand square is given by adjointness and the right-hand square
by applying $U_\alpha$ to the right-hand square of the previous
diagram.  Then there exists a lift $U_\alpha B \rightarrow
U_\alpha X$. Since any fibration of simplicial sets occurs as
$U_\alpha X \rightarrow U_\alpha Y$ for some $X \rightarrow Y$ in
$\Algt$, the map $U_\alpha A \rightarrow U_\alpha B$ is an acyclic
cofibration of simplicial sets and in particular a weak
equivalence.
\end{proof}

We now need a model category structure on the category of homotopy
$\mathcal T$--algebras. However, the category of homotopy $\mathcal
T$--algebras does not have all small limits and colimits (axiom
MC1). Thus, we instead define a model category structure on all
diagrams $\mathcal T \rightarrow \SSets$ in such a way that the
fibrant objects are homotopy $\mathcal T$--algebras.

The following theorem holds for model categories $\mathcal M$
which are left proper and cellular. We will not define these
conditions here, but refer the reader to \cite[13.1.1,
12.1.1]{hirsch} for more details.  It can be shown that $\SSetst$
satisfies both these conditions \cite[13.1.14, 12.5.1]{hirsch}.

\begin{theorem}{\rm\cite[4.1.1]{hirsch}}\label{Loc}\qua
Let $\mathcal M$ be a left proper cellular model category and $P$
a set of morphisms of $\mathcal M$. There is a model category
structure $\mathcal L_P \mathcal M$ on the underlying category of
$\mathcal M$ such that:
\begin{enumerate}[\upshape(1)]
\item The weak equivalences are the $P$--local equivalences.

\item The cofibrations are precisely the cofibrations of $\mathcal
M$.

\item The fibrations are the maps which have the right lifting
property with respect to the maps which are both cofibrations and
$P$--local equivalences.

\item The fibrant objects are the $P$--local objects.
\end{enumerate}
\end{theorem}

To localize the model structure $\SSetst_f$, we first need an
appropriate map. To do so for ordinary algebraic theories,
Badzioch uses free diagrams which are corepresented by the objects
$T_n$ of the theory $\mathcal T$ \cite[2.9]{bad}. In particular
the $n$ projection maps $T_n \rightarrow T_1$ induce maps
\[ \coprod_{i=1}^n \Hom_{\mathcal T}(T_1, -) \rightarrow \Hom_{\mathcal T}(T_n, -). \]
He defines his localization with respect to the set of these maps.
We would like to define similar free diagrams in a multi-sorted
theory.

For each $\alphau^n=<\alpha_1, \ldots ,\alpha_n>$ and $1 \leq i
\leq n$, there exists a projection map $T_{\alphau^n} \rightarrow
T_{\alpha_i}$ inducing a map
\[ \Hom_{\mathcal T}(T_{\alpha_i}, -) \rightarrow \Hom_{\mathcal
T}(T_{\alphau^n}, -). \] Taking the coproduct of all such maps
results in a map
\[ p_{\alphau^n}\co \coprod_{i=1}^n \Hom_{\mathcal
T}(T_{\alpha_i}, -) \rightarrow \Hom_{\mathcal T}(T_{\alphau^n},-).
\]  These maps are the ones which we will use to localize
$\SSetst$.  We define $P$ to be the set of all such maps
$p_{\alphau^n}$ for each $\alphau^n$ and $n \geq 0$.

\begin{prop}
There is a model category structure $\mathcal{LSS}ets^\mathcal T$
on the category $\SSetst$ with weak equivalences the $P$--local
equivalences, cofibrations as in $\SSetst_f$, and fibrations the
maps which have the right lifting property with respect to the
maps which are cofibrations and weak equivalences.
\end{prop}

\begin{proof}
This proposition is a special case of \fullref{Loc}.
\end{proof}

The following propositions are proved by Badzioch for ordinary
theories.  His proofs can be generalized to apply to multi-sorted
theories as well.

\begin{prop}{\rm\cite[5.5]{bad}}\qua
An object $Z$ of $\LSSetst$ is fibrant if and only if it is a
homotopy $\mathcal T$--algebra which is fibrant as an object of
$\SSetst_f$.
\end{prop}

\begin{prop}{\rm\cite[5.6]{bad}}\label{wkequiv}\qua
If $Z$ and $X'$ are homotopy $\mathcal T$--algebras in $\SSetst$
and there is a $P$--local weak equivalence $f\co Z \rightarrow X'$,
then $f$ is also a weak equivalence in $\SSetst_f$, ie, an
objectwise weak equivalence.
\end{prop}

\begin{prop}{\rm\cite[5.8]{bad}}\label{fibrant}\qua
A map $f\co X \rightarrow X'$ is a $P$--local equivalence if and only
if for any $\mathcal T$--algebra $Y$ which is fibrant in
$\SSetst_c$, the induced map of function complexes
\[ f^*\co \Map(X',Y) \rightarrow \Map(X,Y) \]
is a weak equivalence of simplicial sets.
\end{prop}

These results can actually be stated in more generality; they are
really just statements about the fibrant objects in a localized
model category structure (see chapter 3 of \cite{hirsch} for more
details).

Hence, we can consider the category $\LSSetst$ to be our homotopy
$\mathcal T$--algebra model category structure.

\section{Rigidification of algebras over multi-sorted theories}

We are now able to prove the following statement, which is a
stronger version of \fullref{first}:

\begin{theorem} \label{second}
There is a Quillen equivalence of model categories between $\Algt$
and $\LSSetst$.
\end{theorem}

We begin with the necessary definitions.

\begin{definition}\cite[1.3.1]{hovey}\qua
If ${\mathcal C}$ and $\mathcal D$ are model categories, then the
adjoint pair $(F,R, \varphi)$ is a \emph{Quillen pair} if one of
the following equivalent statements is true:
\begin{enumerate}
\item $F$ preserves cofibrations and acyclic cofibrations.

\item $R$ preserves fibrations and acyclic fibrations.
\end{enumerate}
\end{definition}

The following theorem is useful for showing that we have a Quillen
pair of localized model category structures.

\begin{theorem}{\rm\cite[3.3.20]{hirsch}}\label{LocPair}\qua
Let $\mathcal C$ and $\mathcal D$ be left proper, cellular model
categories and let $(F,R, \psi)$ be a Quillen pair between them.
Let $S$ be a set of maps in $\mathcal C$ and $L_S \mathcal C$ the
localization of $\mathcal C$ with respect to $S$.  Then if ${\bf
L}FS$ is the set in $\mathcal D$ obtained by applying the left
derived functor of $F$ to the set $S$ \cite[8.5.11]{hirsch}, then
$(F,R, \psi)$ is also a Quillen pair between the model categories
$L_S\mathcal C$ and $L_{{\bf L}FS}D$.
\end{theorem}

\begin{definition}{\rm\cite[1.3.12]{hovey}}\qua
A Quillen pair is a \emph{Quillen equivalence} if for all
cofibrant $X$ in $\mathcal C$ and fibrant $Y$ in $\mathcal D$, a
map $f\co FX \rightarrow Y$ is a weak equivalence in $\mathcal D$ if
and only if the map $\varphi f\co X \rightarrow RY$ is a weak
equivalence in $\mathcal C$.
\end{definition}

We need to find an adjoint pair of functors between $\Algt$ and
$\LSSetst$ and prove that it is a Quillen equivalence.  Let
\[ J_{\mathcal T}\co  \Algt \rightarrow
\mathcal{SS}ets^\mathcal T \] be the inclusion functor.  We need
to show we have an adjoint functor taking an arbitrary diagram in
$\SSetst$ to a $\mathcal T$--algebra. We first make the following
definition.

\begin{definition}
Let $\mathcal D$ be a small category and $\SSetsd$ the category of
functors $\mathcal D \rightarrow \SSets$. Let $P$ be a set of
morphisms in $\SSetsd$.  An object $Y$ in $\SSetsd$ is
\emph{strictly} $P$--\emph{local} if for every morphism $f\co A
\rightarrow B$ in $P$, the induced map on function complexes
\[ f^*\co  \Map (B,Y) \rightarrow \Map (A,Y) \]
is an isomorphism of simplicial sets. A map $g\co C \rightarrow D$ in
$\SSetsd$ is a \emph{strict} $P$--\emph{local equivalence} if for
every strictly $P$--local object $Y$ in $\SSetsd$, the induced map
\[ g^*\co \Map (D,Y) \rightarrow \Map(C,Y) \]
is an isomorphism of simplicial sets.
\end{definition}

Now, given a category of $\mathcal D$--diagrams in $\SSets$ and the
full subcategory of strictly $P$--local diagrams for some set $P$
of maps, we have the following result.  (Ad\'{a}mek and
Rosick\'{y} also prove this fact \cite[1.38]{ar}, using slightly
different terminology.)

\begin{lemma} \label{adjoint}
Consider two categories, the category of all diagrams $X\co  \mathcal
D \rightarrow \SSets$ and the category of strictly local diagrams
with respect to the set of maps $P= \{f\co A \rightarrow B\}$. Then
the forgetful functor from the category of strictly local diagrams
to the category of all diagrams has a left adjoint.
\end{lemma}

\begin{proof}
Without loss of generality, assume that we have just one map $f$
in $P$; otherwise replace $f$ by $\coprod_\alpha f_\alpha$. Given
an arbitrary diagram $X$, we would like to construct a strictly
local diagram from $X$. So, suppose that $X$ is not strictly
local, ie, the map
\[ f^*\co \Map(B,X) \rightarrow \Map(A,X) \]
is not an isomorphism.  To ensure that $f^*$ is surjective, we
obtain an object $X'$ as the pushout in the following diagram:
\[ \xymatrix{\coprod_{n \geq 0} \coprod_{A \times \Delta [n]
\rightarrow X} A \times \Delta[n] \ar[r] \ar[d] & X \ar[d] \\
\coprod_{n \geq 0} \coprod_{A \times \Delta [n] \rightarrow X} B
\times \Delta [n] \ar[r] & X'}
\] where each coproduct is taken over all maps $A \times \Delta
[n] \rightarrow X$ for each $n \geq 0$. Then, to ensure that $f$
is injective, we obtain $X''$ by taking the pushout
\[ \xymatrix{\coprod_{n \geq 0} \coprod (B \coprod_A B) \times
\Delta [n] \ar[r] \ar[d] & X' \ar[d] \\
\coprod_{n \geq 0} \coprod B \times \Delta [n] \ar[r] & X''} \]
again where the second coproduct is over all maps $(B \coprod_A B)
\times \Delta [n] \rightarrow X'$, and where the map
\[ B \coprod_A B \rightarrow B \]
is the fold map.
\vspace{2pt}

In the construction of $X'$, for any strictly local object $Y$ we
obtain a pullback diagram
\[ \xymatrix{\Map(X',Y) \ar[r] \ar[d]^\cong & \Map(\coprod B,Y)
\ar[d]^\cong \\
\Map(X,Y) \ar[r] & \Map(\coprod A,Y)} \] showing that the map $X
\rightarrow X'$ is a strict local equivalence since $f\co A
\rightarrow B$ is.
\vspace{2pt}

In the construction of $X''$, we obtain a similar diagram, but it
takes more work to show that the map $X' \rightarrow X''$ is a
strict local equivalence.  We first obtain the pullback diagram
\[ \xymatrix{ \Map(X'',Y) \ar[r] \ar[d] & \Map(\coprod B,Y) \ar[d] \\
\Map(X',Y) \ar[r] & \Map(\coprod(B \coprod_A B),Y)} \] Since it is
a pullback diagram then suffices to show that the right hand
vertical arrow is an isomorphism.

Recall that the object $B \coprod_AB$ is defined as the pushout in
the diagram
\[ \xymatrix{A \ar[r] \ar[d] & B \ar[d] \\
B \ar[r] & B \coprod_AB} \] which enables us to look at the
pullback diagram
\[ \xymatrix{ \Map(B \coprod_AB,Y) \ar[r] \ar[d] & \Map(B,Y) \ar[d]^\cong \\
\Map(B,Y) \ar[r] & \Map(A,Y). } \] Hence the map
\[ B \rightarrow B \coprod_A B \]
is a strict local equivalence.  But, this map fits into a
composite
\[ \xymatrix@1{B \ar[r] \ar@/_1pc/[rr]_{\rm id} & B \coprod_A B \ar[r]
& B} \]  Since the identity map is a strict local equivalence, it
follows that the map
\[ B \coprod_A B \rightarrow B \]
is a strict local equivalence, since it can be shown that the
strictly local equivalences satisfy model category axiom MC2.

Therefore, we obtain a composite map $X \rightarrow X''$ which is
a strict local equivalence. However, we still do not know that the
map
\[ \Map(B,X'') \rightarrow \Map(A,X'') \]
is an isomorphism.  So, we repeat this process, taking a (possibly
transfinite) colimit to obtain a strictly local object $\widetilde
X$ such that there is a local equivalence $X \rightarrow
\widetilde X$.

It suffices to show that the functor which takes a diagram $X$ to
the local diagram $\widetilde X$ is left adjoint to the forgetful
functor.  So if $J$ is the forgetful functor from the category of
strictly local diagrams to the category of all diagrams and $K$ is
the functor we have just defined, we claim that
\[ \Map(X, JY) \cong \Map(KX, Y) \]
for any diagram $X$ and strictly local diagram $Y$.  But, proving
this statement is equivalent to showing that
\[ \Map(X,Y) \cong \Map(\widetilde X,Y) \]
which was shown above for each step, and it still holds for the
colimit.  In particular, the map $X \rightarrow \widetilde X=KX$
and the identity $Y=JY$ induce natural isomorphisms
\[ \Map(KX,Y) \rightarrow \Map (X,Y) \rightarrow \Map(X,JY), \]
and the restriction of this composite to the 0--simplices of each
object,
\[ \Hom(KX,Y) \rightarrow \Hom(X,JY) \] is exactly the
isomorphism we need to show that $K$ is left adjoint to $J$.
\end{proof}

To apply this lemma to our situation, we first need to verify that
$\Algt$ is precisely the category of strictly local diagrams in
$\SSetst$ with respect to the set of maps $P$, defined in the last
section, to obtain the model category structure for homotopy
$\mathcal T$--algebras.

To do so, we will need the following homotopical version of the
Yoneda Lemma.

\begin{lemma} \label{yoneda}
Let $T_{\alphau}$ be an object in a multi-sorted algebraic theory
$\mathcal T$ and $A$ a strict $\mathcal T$--algebra.  There is an
isomorphism of simplicial sets
\[ \Map_{\SSetst}(\Hom_{\mathcal T}(T_\alphau, -), A) \cong A(T_\alphau). \]
\end{lemma}

\begin{proof}
Since $A$ is a simplicial set-valued functor, we can regard it as
a simplicial diagram of set-valued functors
\[ A(-)_0 \Leftarrow A(-)_1 \Lleftarrow A(-)_2 \cdots \]
which further induces a simplicial diagram
\[ \Hom_{\SSetst}(\Hom_{\mathcal T} (T_\alphau, -), A(-)_0)
\Leftarrow \Hom_{\SSetst}(\Hom_{\mathcal T}(T_\alphau, -), A(-)_1)
\Lleftarrow \cdots. \]  Using the classical Yoneda Lemma
\cite[11.5.8]{hirsch}, we have a natural isomorphism at each level
\[ \Hom_{\Sets^\mathcal T}(\Hom_{\mathcal T}(T_\alphau, -), A(-)_n) \cong
A(T_\alphau)_n, \] where $\Sets^\mathcal T$ denotes the category
of functors $\mathcal T \rightarrow \Sets$.

Now, regarding sets as constant simplicial sets as necessary,
notice that there are natural isomorphisms
\[ \begin{aligned}
\Hom_{\Sets^\mathcal T}&(\Hom_{\mathcal T}(T_\alphau, -), A(-)_n) \\ &
\cong \Hom_{\SSetst}(\Hom_{\mathcal T}(T_\alphau,-),
\Hom_{\SSetst}(\Delta
[n], A(-))) \\
& \cong \Hom_{\SSetst}(\Hom_{\mathcal T}(T_\alphau, -) \times \Delta
[n], A(-)) \\
& \cong \Map_{\SSetst}(\Hom_{\mathcal T}(T_\alphau,-), A)_n.
\end{aligned} \]
Since all the simplicial maps above are natural, we obtain a
natural simplicial functor
\[ \Map_{\SSetst}(\Hom_{\mathcal T}(T_\alphau,-),A) \rightarrow
A(T_\alphau) \] which is an isomorphism.
\end{proof}

Using this lemma, we are able to prove the following.

\begin{lemma} \label{strictlylocal}
A diagram $A\co \mathcal T \rightarrow \SSets$ is a strict $\mathcal
T$--algebra if and only if $A$ is strictly local with respect to
the maps
\[ p_{\alphau^n}\co  \coprod_{i=1}^n \Hom_{\mathcal T}(T_{\alpha_i},-) \rightarrow \Hom_{\mathcal
T}(T_{\alphau},-). \]
\end{lemma}

\begin{proof}
A diagram $A$ is a strict $\mathcal T$--algebra if and only if for
each $\alphau^n=<\alpha_1, \ldots ,\alpha_n>$ there is a natural
isomorphism
\[ \prod_{i=1}^n A(T_{\alpha_i}) \cong  A(T_\alphau) \] induced by
the projection maps in $\mathcal T$.  Using \fullref{yoneda},
this statement is equivalent to having an isomorphism
\[ \begin{aligned}
\Map_{\SSetst}(\Hom_{\mathcal T}(T_\alphau,-), A) & \cong
\prod_{i=1}^n \Map_{\SSetst}(\Hom_{\mathcal T}(T_{\alpha_i},-), A)
\\
& \cong \Map_{\SSetst}(\coprod_{i=1}^n \Hom_{\mathcal
T}(T_{\alpha_i}, -), A)
\end{aligned}. \]  Since all the isomorphisms in sight are induced
by projections, it follows that this statement is equivalent to
having $A$ strictly local with respect to all the maps
$p_{\alphau^n}$.
\end{proof}

In particular, a map $f\co X \rightarrow X'$ is a strict $P$--local
equivalence if and only if for every $A$ in $\Algt$ (regarded as
an object in $\SSetst$ via the map $J_{\mathcal T}$) the induced map
\[ \Map_{\SSetst}(X',A) \rightarrow \Map_{\SSetst}(X,A) \] is an isomorphism of
simplicial sets.

Applying \fullref{adjoint} to the functor $J_{\mathcal T}$, we
obtain its left adjoint functor
\[ K_{\mathcal T}\co \SSetst \rightarrow \Algt. \]

\begin{prop}
The adjoint pair of functors
\[ \xymatrix@1{K_{\mathcal T}\co  \SSetst \ar@<.5ex>[r] & \Algt \co J_{\mathcal T}. \ar@<.5ex>[l]} \]
is a Quillen pair.
\end{prop}

\begin{proof}
Using \fullref{strictlylocal}, we can regard $\Algt$ as a
subcategory of $\SSetst$ via the map $J_{\mathcal T}$. Since in both
cases, the fibrations and weak equivalences are defined
objectwise, $J_{\mathcal T}$ preserves fibrations and acyclic
fibrations.
\end{proof}

\begin{lemma} \label{ktk}
Each map $K_{\mathcal T}(p_{\alphau^n})$ is a weak equivalence in
$\Algt$.
\end{lemma}

\begin{proof}
First, we note that the functor $\Hom_{\mathcal T}(T_\alphau,-)$ is
a strict $\mathcal T$--algebra and that $J_{\mathcal T} A=A$ for any
strict $\mathcal T$--algebra $A$, again regarding $J_{\mathcal T}$ as
an inclusion functor. Then, for each map $p_{\alphau^n}$ we have
the following composite isomorphism:
\[ \begin{aligned}
\Map_{\Algt}(K_{\mathcal T}(\Hom_{\mathcal T}(T_\alphau,-)),A) & \cong
\Map_{\SSetst}(\Hom_{\mathcal T}(T_\alphau,-),A) \\
& \cong A(T_\alphau) \\
& \cong \prod_{i=1}^n A(T_{\alpha_i}) \\
& \cong \prod_{i=1}^n \Map_{\SSetst}(\Hom_{\mathcal
T}(T_{\alpha_i},-), A) \\
& \cong \Map_{\SSetst}(\coprod_{i=1}^n \Hom_{\mathcal
T}(T_{\alpha_i},-), A) \\
& \cong \Map_{\Algt}(K_{\mathcal T}(\coprod_{i=1}^n \Hom_{\mathcal
T}(T_{\alpha_i},-)), A).
\end{aligned} \]  Since all the isomorphisms are naturally induced by
the map $p_{\alphau^n}$ and adjoints, it follows that $K_{\mathcal
T}$ is a strict local equivalence, or a weak equivalence in
$\Algt$.
\end{proof}

Now, we need to show that the same adjoint pair is still a Quillen
pair when we replace the model structure $\SSetst$ with the model
structure $\LSSetst$.

\begin{prop}
The adjoint pair
\[ \xymatrix@1{K_{\mathcal T}\co  \LSSetst \ar@<.5ex>[r] & \Algt
\,\,\colon\! J_{\mathcal T} \ar@<.5ex>[l]} \] is a Quillen pair.
\end{prop}

\begin{proof}
Consider again the set of maps
\[ P= \{p_{\alphau^n}\co \coprod_i \Hom_{\mathcal T}(T_{\alpha_i},-) \rightarrow \Hom_{\mathcal T}(T_{\alphau},-)\}. \]
Notice in particular that the objects involved in these maps are
free diagrams and therefore cofibrant in $\SSetst_f$.  The model
category structure $\LSSetst$ is obtained by localizing with
respect to these maps. Then using \fullref{ktk}, we have that
each map $K_{\mathcal T} (p_{\alphau^n})$ is a weak equivalence in
$\Algt$. Hence, it follows from \fullref{LocPair} that the
pair of adjoints forms a Quillen pair even after the localization
on $\SSetst_f$.
\end{proof}

Before stating the main theorem, that the above Quillen pair is
actually a Quillen equivalence, we first need a lemma.  Badzioch's
proof \cite[6.5]{bad} for ordinary theories generalizes for our
case of multi-sorted theories, but we give a slightly different
proof here.

\begin{lemma}
If $X$ is cofibrant in $\LSSetst$, then the unit map $\eta \co X
\rightarrow K_{\mathcal T} X = J_{\mathcal T} K_{\mathcal T} X$ is a
weak equivalence in $\LSSetst$.
\end{lemma}

\begin{proof}
Case 1: The cofibrant object $X$ is a free diagram, so it can be
written as
\[\coprod_{\alphau} \Hom_{\mathcal T}(T_{\alphau^n},-). \]  The proof
for such an object is then similar to the argument in the proof of
\fullref{ktk}.

Case 2: Let $X$ be any cofibrant diagram.  Then $X \simeq
\hocolim_{\Deltaop} X_i$ where each $X_i$ is a free diagram. Using
\fullref{fibrant}, it then suffices to show that
$\Map(K_{\mathcal T} X,Y) \simeq \Map(X,Y)$ for any $\mathcal
T$--algebra $Y$ which is fibrant in $\SSetst_{cof}$. Using case 1,
we have the following:
\[ \begin{aligned}
\Map(X,Y) & \simeq \Map(\hocolim_{\Deltaop} X_i,Y) \\
& \simeq \holim_{\bf \Delta} \Map(X_i,Y) \\
& \simeq \holim_{\bf \Delta} \Map(K_{\mathcal T} X_i,Y) \\
& \simeq \Map(\hocolim_{\Deltaop}K_{\mathcal T} X_i,Y) \\
& \simeq \Map(K_{\mathcal T} X,Y).
\end{aligned} \]  Notice in particular that this weak equivalence
is induced by the map $\eta$. The lemma follows.
\end{proof}

Now, the proof of the main theorem follows from this lemma exactly
as it does for ordinary theories in \cite[6.4]{bad}.

\begin{theorem}
The Quillen pair of functors
\[ \xymatrix@1{K_{\mathcal T}\co  \LSSetst \ar@<.5ex>[r] & \Algt
\,\,\colon\! J_{\mathcal T}. \ar@<.5ex>[l]} \] is a Quillen equivalence.
\end{theorem}

\begin{proof}
Let $X$ be a cofibrant object in $\LSSetst$, $A$ a fibrant object
in $\Algt$, and $f\co X \rightarrow A = J_{\mathcal T} A$ a map in
$\LSSetst$.  We need to show that $f$ is a $P$--local equivalence
if and only if its adjoint map $g\co K_{\mathcal T} X \rightarrow A$ is
a weak equivalence in $\Algt$.  There is a commutative diagram
\[ \xymatrix{X \ar[r]^-\eta \ar[d]_f & K_{\mathcal T} X \ar[d]^g  \\
J_{\mathcal T}A \ar[r]^-= & A} \]
First assume that $f$ is a $P$--local equivalence.  Then $g$ must
also be a $P$--local equivalence since $\eta$ is, by the previous
lemma.  However, $g$ is a map in $\Algt$, and so it is an
objectwise weak equivalence, or a weak equivalence in $\Algt$.

Conversely, suppose that $g$ is a weak equivalence in $\Algt$.
Then it is a $P$--local equivalence.  Hence, $f = g\circ \eta$ is
also a $P$--local equivalence.
\end{proof}

Hence, we have a Quillen equivalence of model categories between
strict $\mathcal T$--algebras and homotopy $\mathcal T$--algebras.

\bibliographystyle{gtart}
\bibliography{link}

\end{document}